\documentclass[final,1p,times,authoryear]{elsarticle}
\usepackage[english]{babel}

\usepackage[T1]{fontenc}

\usepackage{graphicx}

\usepackage{amsthm}
\usepackage{color}

\usepackage{caption}
\usepackage{tabularx}

\usepackage[]{hyperref}

\usepackage{latexsym}
\usepackage{amsmath}
\usepackage{graphics}
\usepackage{amssymb}
\usepackage{color}
\usepackage{subcaption}
\newtheorem{theorem}{Theorem}

\newtheorem*{theorem*}{Theorem}

\newtheorem{lemma}[theorem]{Lemma}

\newtheorem{remark}[theorem]{Remark}

\newtheorem{proposition}[theorem]{Proposition}

\newtheorem{definition}[theorem]{Definition}

\newtheorem*{problem1*}{Problem A} 
\newtheorem*{problem2*}{Problem B}

\DeclareMathOperator{\inter}{int} 
\DeclareMathOperator{\diver}{div}

\DeclareMathOperator{\cof}{Cof}

\newcommand{\beq}{\begin{equation}}
\newcommand{\eeq}{\end{equation}}

\newcommand{\la}{{\lambda}}

\newcommand{\La}{{\Lambda}}

\newcommand{\Om}{{\Omega}}

\newcommand{\D}{\nabla}

\newcommand{\sio}{\sigma_i^0}

\newcommand{\s}{\sigma}

\newcommand{\ph}{\varphi}

\newcommand{\po}{{\partial\Omega}}

\newcommand{\ny}{{\nabla y}}

\newcommand{\bpr }{\begin{proof}}
\newcommand{\epr}{\end{proof}}

\newcommand{\hw}{{\widehat W}}

\newcommand{\hs}{{\widehat S}}

\newcommand{\cD}{{\mathcal{D}}}

\newcommand{\cS}{{\mathcal{L}}}

\newcommand{\cC}{{\mathcal{C}}}

\newcommand{\cA}{{\mathcal{A}}}

\newcommand{\re}{\mathbb{R}}

\newcommand{\ren}{{\mathbb{R}^N}}

\newcommand{\renp}{{\mathbb{R}_+^N}}

\newcommand{\rt}{{\mathbb{R}^2}}
\newcommand{\lp}{{\mathrm{Lin}_+}}

\newcommand{\ort}{{\mathrm{Orth}_+}}

\newcommand{\sym}{{\mathrm{Sym}}}

\biboptions{round}

\numberwithin{equation}{section}

\numberwithin{theorem}{section}

\begin{document}

\begin{frontmatter}



\title{Uniqueness for Some Mixed Problems of Nonlinear Elastostatics}


\author[]{Phoebus Rosakis}
\ead{rosakis@uoc.gr}

%\address[label1]
 \address{Department of Mathematics and Applied Mathematics, University of Crete, Heraklion 70013 Crete, Greece}
 \address{Institute for Applied and Computational Mathematics, \\ Foundation of Research and Technology-Hellas, Heraklion 70013, Crete, Greece}\bigskip
\address{Center of Mediocrity for Antiquated Materials, Mitato, Crete 72300, Greece}
\address{{.}}
\address{{\normalsize Dedicated to Rohan Abeyaratne on the occasion of his 70th birthday}}

\begin{abstract}
We show that certain  mixed displacement/traction problems (including live pressure tractions) of nonlinear elastostatics that are solved by a homogeneous deformation, admit no other classical equilibrium solution under suitable constitutive inequalities and domain boundary restrictions.  This extends a well known theorem of  Knops and Stuart on the pure displacement problem.
\end{abstract}

%\begin{keyword}
%Science \sep Publication \sep Complicated


%\end{keyword}

\end{frontmatter}

% _____________________________________________________________________________

%\tableofcontents
%\clearpage

% _______________________________________________________________

\section{Introduction}

\cite{ks} have established global uniqueness of classical equilibrium solutions of the \emph{homogeneous} displacement problem of nonlinear elastostatics (where boundary displacements are prescribed consistent with a given \emph{homogeneous} deformation).  They require the reference region $\Om\subset\ren $ to be star shaped, the stored energy function $W$ to be (globally) rank-one convex, and  strictly quasiconvex at the prescribed homogeneous deformation (see Definition \ref{d1}). This seems to be the first global uniqueness result for nonlinear elastostatics under physically acceptable constitutive restrictions, stated here in Theorem \ref{t3}. For various uniqueness theorems for nonlinear elastostatics, see the works cited by \cite{ks} and \cite{sis}. 

Here we provide uniqueness results analogous to that of \cite{ks} for the homogeneous mixed problem, where displacements and tractions (consistent with a given homogeneous deformation) are prescribed on two complementary subsets $\cD$ and $\cS$ of the boundary $\po$, respectively. 
After some preliminaries in Section \ref{s2}, we study in Section \ref{s3} the mixed displacement / dead-load traction problem.
\begin{problem1*} 
Given $F_0\in\lp$, find a classical equilibrium $y:\Om\to\ren$ subject to the boundary conditions of displacement
$$ y(x)=F_0x \quad \forall x \in \cD $$
and dead load traction
 $$ S(\ny)n=S(F_0)n \quad \hbox{on } \cS $$
where $n$ is the unit outward normal on $\po$.
\end{problem1*}
Here $\lp$ contains all tensors with positive determinant, $S(F)=DW(F)$ is the Piola stress, a classical equilibrium is a classically smooth deformation $y$ such that $S(\ny)$ is divergence free; see definition \ref{d1}. Uniqueness is addressed in Propositions \ref{l4}, \ref{l5}, \ref{l7}. 

The presence of additional traction boundary conditions is dealt with using similar tools and methods as in \cite{ks}. In particular, Green's Identity (Lemma \ref{l2}) plays a central role, but we require different restrictions on boundary geometry and more stringent constitutive inequalities. 

While $\Om$ need not be star shaped, the boundary $\po$ is required to satisfy a \emph{partition condition} (Definition \ref{coo}) that restricts how $\po$  splits into complementary subsets $\cD$ and $\cS$, where displacements and tractions are prescribed. The Partition Condition is the requirement that $x\cdot n\ge 0$ for $x\in \cD$ and $x\cdot n\le 0$ on $\cS$, with $n$ the outward unit normal at $x\in\po$. Various examples are shown in Fig.~1 and discussed in Section \ref{61}.

The stored energy $W$ is required to be (sometimes strictly) rank-one convex, but quasiconvexity at the gradient $F_0$ of the homogeneous deformation is apparently insufficient for uniqueness, due to the presence of traction conditions on part of $\po$.

Instead, in most of our results we demand that $F_0$ be a point of convexity of $W$. Although global convexity of $W:\lp\to\re$  is inappropriate for nonlinear elasticity (\cite{b}), nonconvex functions can have many points of convexity. In Section \ref{62} we show a stored energy function of compressible Neo-Hookean type, in two dimensions, where all deformation gradients $F$ with $\det F>1$ are points of strict convexity of $W$, while if $\det F<1$, $F$ is not even a point of convexity of $W$.  In particular, states of uniaxial tension of a slightly tapered slender column as in Fig.~1g are covered by our uniqueness result, Proposition \ref{l7}, but uniaxial compression is not, in accordance with the expectation of  nonuniqueness due to buckling. 

We specialize $W$ to be isotropic in Section \ref{s4}. This allows somewhat weaker convexity hypotheses at the homogeneous state, but these require a much more delicate proof. See Proposition \ref{l7}. As a result, natural states that are not points of strict convexity due to rotational invariance, uniaxial tension in three dimensions, also certain stored energies with multiple wells that are nonetheless polyconvex (with  wells that are  not rank-one connected) satisfy our hypotheses.

In Section \ref{s5}, we investigate  a mixed displacement/Cauchy pressure problem, an example of live loads. The homogenenous deformation is now restricted to be a dilatation: $F_0=J_0^{1/N}I$ with prescribed $J_0=\det F_0$.
\begin{problem2*}
Given $J_0>0$, find a classical equilibrium $y$ subject to the boundary conditions of displacement
$$ y(x)= J_0^{1/N} x \quad \forall x \in \cD $$
and traction
$$  S(\ny)n=\ph'(J_0)\cof (\ny)n \quad \hbox{on } \cS $$
where $n$ is the unit outward normal on $\po$.
\end{problem2*} 
Here $\ph(J):=W(J^{1/N}I)$ for $J>0$ is the restriction of the energy to dilatations.
The traction boundary condition on $\cS$ is equivalent to prescribing the  Cauchy traction, in this case a hydrostatic pressure or tension corresponding to the given dilatation, on the deformed boundary $y(\cS)$. This allows considerably weaker constitutive restrictions; $J_0$ must be a point of convexity of  $\ph$. This is actually weaker than requiring $F_0$ to be a point of convexity of $W$ (Lemma \ref{l6} due to \cite{s}), in case the Cauchy stress is a pressure (negative principal Cauchy stress). See Proposition \ref{l9} for uniqueness of solutions to Problem B.

Section 6 contains a discussion of the various hypotheses we require for uniqueness and examples of regions $\Om$, stored energies $W$ and homogeneous states $F_0$ that satisfy or violate them, in comparison with reasonable expectations for uniqueness. Specifically, we present an explicit example where uniqueness fails when the Partition Condition is omitted from our hypotheses.

We note that \cite{sis} have shown uniqueness of equilibria for the mixed problem with dead loads, with virtually no restrictions on $\cD$ and $\cS$, while assuming uniform polyconvexity of $W$. On the other hand, their equilibria are restricted to satisfy additional conditions in the form of certain pointwise bounds. They also have various references on uniqueness. Under our (different) hypotheses, we show by example that our uniqueness results are false  when the Partition Condition is omitted from our hypotheses; see Section \ref{63}.

\section{Preliminaries}\label{s2} \noindent Let $\Om\subset\ren$ be a compact region with nonempty interior $\inter\Om$ and  $\po=\cD\cap\cS $, where $\cD$ and $\cS $ are disjoint, piecewise smooth (hyper)surfaces and $\cD$ is  nonempty. The stored energy function is $W\in C^2(\lp,\re)$. Its derivative is the Piola stress tensor (more pedantically called the \emph{first Piola-Kirchhoff stress tensor})
$$ S(F)=DW(F) \quad \forall F\in \lp.$$
\begin{definition}\label{d1}(i) The stored energy function $ W \colon \lp\to \re $ is rank-one convex provided  
\beq W(A)-W(F)\ge S(F)\cdot (A-F) \label{CON}\eeq
holds $\forall \;A$, $F \in\lp$ with $ A-F=a\otimes b$, $a$, $b\in\ren$. \hfill\break
(ii) $F_0\in\lp$ is a point of convexity of $W$ if \eqref{CON} holds with $F=F_0$ and $\forall A\in\lp$.
\hfill\break 
(iii) $W$ is quasiconvex at $F_0\in\lp$ provided
$$\int_\Om W(F_0+\D v(x))dx\ge |\Om|W(F_0)\;\; \forall v\in C^1_0(\Om,\ren)$$
(iv)A classical equilibrium is a mapping $y\in C^2(\inter\Om,\ren)\cap C^1(\bar\Om,\ren)$ with $\D y\in\lp$ on $\Om$, satisfying 
\beq \diver   S(\D y)=0 \hbox{ on }\inter\Om.\label{cleq}\eeq

\end{definition}

\begin{lemma}\label{l2} (Green's Identity) If $y$ is a classical equilibrium, then
\beq N\int_\Om W(\D y)dx=\int_\po \Bigl\{(x\cdot n)W(\D y) +S(\D y)n\cdot[y-(\D y)x]\Bigr\}dS.\label{gi}\eeq
\end{lemma}
\bpr Assume $y$ is a classical equilibrium, so that the Eshelby tensor \beq P:=W(\D y)I-(\D y)^TS(\D y)\label{et}\eeq is divergence free, which implies that 
$$\diver (P^Tx)=\rm{tr} P=NW-S\cdot\ny.
$$
Integrate this over $\Om$, use the Divergence Theorem
and note that
$$\int_\Om S\cdot\D ydx=\int_\po Sn\cdot ydS$$
in view of \eqref{cleq}, to obtain
$$N\int_\Om W(\D y)dx=\int_\po (Pn\cdot x+ Sn\cdot y)dS$$
from which \eqref{gi} follows after using \eqref{et}.\epr 

\section{The Mixed Displacement/Dead-Load Problem}\label{s3}
\noindent We consider the following elastostatic problem with mixed displacement and dead-load traction boundary conditions consistent with a homogeneous deformation:
\begin{problem1*} 
Given $F_0\in\lp$, find a classical equilibrium $y$ subject to the boundary conditions of displacement
\beq\label{bcu} y(x)=F_0x \quad \forall x \in \cD \eeq
and dead load traction
\beq \label{bcs} S(\ny)n=S(F_0)n \quad \hbox{on } \cS \eeq
where $n$ is the unit outward normal on $\po$.
\end{problem1*}
\noindent Trivially, $y(x)=F_0x$ for $x\in\bar\Om$ is a classical equilibrium solution of this problem.
In case of the pure displacement problem, \citet{ks} have established
\begin{theorem}\label{t3} [\cite{ks}] Suppose that $\cS=\emptyset$, and (i) $\Om$ is star shaped, (ii) $W$ is rank-one convex, (iii) $W$ is  strictly quasiconvex at $F_0$. Then the only classical equilibrium solution of Problem1 is $y=F_0x$ on $\Om$. \end{theorem} 
This makes use of Green's Identity, Lemma \ref{l2} above. Here we extend this result to some mixed problems. Specifically we adopt
\begin{definition}\label{coo} We say that $\cD$ and $\cS$ satisfy the Partition Condition (with respect to the origin),  provided
\beq\label{coopd} x\cdot n(x)\ge0 \quad \forall x\in\cD\eeq
and
\beq \label{coops} x\cdot n(x)\le0 \quad  \forall x\in\cS,\eeq
where $n$
is the unit outward normal on $\po$.
\end{definition}
\begin{remark}
Whether or not $\cD$ and $\cS$ satisfy the Partition Condition depends on the choice of the origin. See  Fig.~1  for various examples where the Partition Condition holds, and Section  \ref{61} for a discussion of this condition.\end{remark}
\begin{figure}
  \centering
\includegraphics[width=0.8\textwidth]{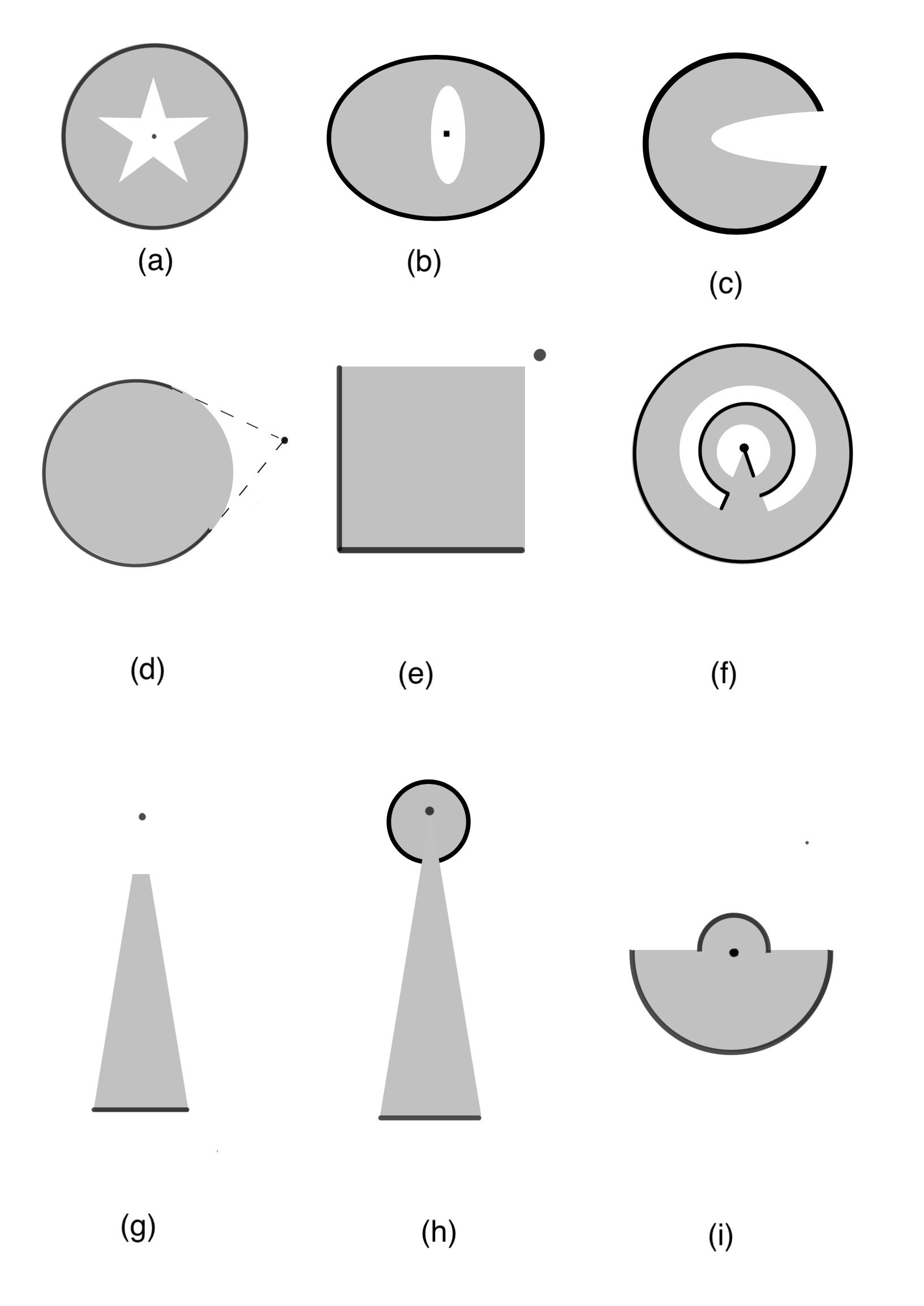}
  \caption{Examples of Partition Condition satisfaction. Here $\Om$ is shown grey. The heavy black portion of the boundary is $\cD$, the remaining part of the boundary between grey and white is $\cS$. The black dot is a possible choice of the origin so that the Partition Condition inequalities in Definition \ref{coo} hold. In all cases except (b,c,d,e), a  revolution about a vertical axis through the origin produces a region in $\re^3$ with the Partition Condition in force. In (f,g,h) the cone boundary points can belong to either $\cD$ or $\cS$. In (c) a solid cube in $\re^3$ with $\cD$ consisting of 3 faces that share a vertex is also admissible.  The tapered column in (g) can be made arbitrarily slender for a fixed height and the taper can be made as small as desired.}
  \label{Fig1}
\end{figure} 

Using essentially the  approach of \cite{ks} we first show

\begin{lemma}\label{l3} Suppose that  (i) $\cD$ and $\cS$ satisfy the Partition Condition (Definition \ref{coo}), (ii) $W$ is rank-one convex, (iii) $W(F)\ge W(F_0)\quad \forall F\in\lp$. Then any  classical equilibrium solution $y$ of Problem A satisfies 
\beq\label{e}W\bigl(\nabla y(x)\bigr)=W(F_0) \quad \forall x\in\Om. \eeq
Moreover if $W$ is strictly rank-one convex, then
\beq\label{ff0}\ny=F_0  \quad\hbox{ on } \cD\subset\po.\eeq
\end{lemma}
\bpr Hypothesis (iii) and the smoothness of $W$ imply that $S(F_0)=0$. Hence, by \eqref{bcs},
\beq\label{nuls} S(\ny)n=0 \hbox{ on }\cS.\eeq
Also, \eqref{bcu} ensures that
\beq\label{r1D} F:=\ny=F_0+a\otimes n  \hbox{ on }\cD.\eeq
Write Green's Identity \eqref{gi} for the two equilibria $y$ and $y_0=F_0x$ and subtract, to obtain, writing $F$ for $\ny$
\beq\label{del}\Delta:=N\int_\Om [W(F)-W(F_0)]dx=I_\cD+I_\cS\eeq
where (using $S(F_0)=0$),
\beq\label{id}I_\cD=\int_\cD \Bigl\{(x\cdot n)\bigl[W(F)-W(F_0)\bigr] +S(F)n\cdot(y-Fx)\Bigr\}dS\eeq
and in view of \eqref{nuls},
\beq\label{is} I_\cS=\int_\cS (x\cdot n)\bigl[W(F)-W(F_0)\bigr]dS\eeq
Appeal to \eqref{bcu}, \eqref{r1D}, to show that
$$ S(F)n\cdot(y-Fx) = S(F)n\cdot(F_0-F)x=$$
$$-S(F)n\cdot[(a\otimes n)x] =-(x\cdot n)S(F)\cdot(a\otimes n)= (x\cdot n)S(F)\cdot(F_0-F)$$
Use this to rewrite $I_\cD$ in \eqref{id} as 
\beq\label{sr1}I_\cD=-\int_\cD (x\cdot n)\bigl[ W(F_0)-W(F)- S(F)\cdot (F_0-F)\bigr]dS.\eeq
In view of Hypothesis (i) (Eq.\eqref{coopd}) and Hypothesis (ii), the integrand is nonnegative on $\cD$, ensuring that $I_\cD\le0$. The above argument is borrowed from \cite{ks}. 

Next, Hypotheses (i)  (Eq.\eqref{coops}), (iii) and  \eqref{is} imply that $I_\cS\le0$.
As a result $\Delta\le0$ in \eqref{del}. Now (iii) ensures that $\Delta\ge0$, thus $\Delta=0$, which proves the result in view of (iii). 

To show \eqref{ff0}, strict rank-one convexity and \eqref{sr1} furnish $I_\cD<0$ if \eqref{ff0} is false, which together with $I_\cS\le0$ contradicts $\Delta=0$.
\epr
In the above, we have made no invariance assumptions on $W$, not even frame indfference, but merely required $F_0$ to minimize it. Next we assume that W is mimimized  on the rotation group $\ort$ (only) consistently with frame indifference, and prove uniqueness under boundary conditions consistent with the natural state, that is, a global minimum of $W$.

\begin{proposition}\label{l4} Suppose that $F_0=1$ and that (i) $\cD$ and $\cS$ satisfy the Partition Condition (Definition \ref{coo}), (ii) $W$ is rank-one convex, (iii) 
\beq\label{rmin}W(F)>W(R)=W(I) \quad \forall F\in\lp\backslash\ort, \quad \forall R\in\ort.\eeq 
Then the only  classical equilibrium solution $y$ of Problem A is $y=x$ on $\Om$
\end{proposition}
\bpr Lemma \ref{l3} applies and \eqref{e} combined with \eqref{rmin}  dictate that $\ny\in\ort$ on $\Om$. Compatibility then demands that $\ny(x)=R=$const. $\in\ort$ $\forall x \in \Om$. 
But then \eqref{r1D} with $F_0=I$ ensures that $R=I$
as no other rotations are rank-one connected to the identity tensor $I$ (\cite{bj}).\epr
We turn to more general choices for the homogeneous state $F_0$ where uniqueness can still be shown. 

\begin{proposition}\label{l5} Suppose that (i) $\cD$ and $\cS$ satisfy the Partition Condition (Definition \ref{coo}), (ii) $W$ is rank-one convex, (iii) $F_0$ is a point of strict convexity of $W$: 
$$W(F)>W(F_0)+S(F_0)\cdot(F-F_0) \quad \forall F\in\lp\backslash\{F_0\}$$
Then the only  classical equilibrium solution $y$ of Problem A is $y=F_0x$ on $\Om$.
\end{proposition}
\bpr Define
\beq\label{hw}\hw (A):= W(A)-W(F_0)-S(F_0)\cdot(A-F_0) \quad \forall A\in\lp .\eeq
Then by (iii),
\beq\label{hmin} \hw(F)>\hw(F_0) \hbox{ for } F\neq F_0\eeq
while
\beq\label{hs}\hs(F):=D\hw(F)=S(F)-S(F_0)\quad\forall F\in\lp.\eeq
Let $y$ be a classical equilibrium solution of Problem A. Then it also a classical equilibrium solution after replacing $W$ with $\hw$ and $S$ by $\hs$ in \eqref{cleq} and Problem A, while Lemma \ref{l3} applies because $F_0$ mninimizes $\hw$. But now \eqref{hmin} together with \eqref{e} ensure that $\ny=F_0$ on $\Om$. The result follows from \eqref{bcu}.\epr

\begin{remark} While global convexity of $W$ is unacceptable as a constitutive inequality in nonlinear elasticity (\cite{b}), in general W has many points of convexity, where it coincides with its convex hull; see Lemma \ref{l6}. For isotropic materials we address  this in  Section \ref{s4} and discuss it in Section \ref{62} where we also provide examples.
\end{remark}

\section{Isotropic Materials}\label{s4}

Let $\renp$
consist of all vectors in $\ren$ with positive entries. The isotropy of $W$ is equivalent to the existence of a symmetric $\Phi:\renp\to\re$ such that 
\beq\label{Phi}W(F)=\Phi(\la(F)) \quad\forall F\in\lp.\eeq where $\la(F)=\bigl(\la_1(F),\ldots,\la_N(F)\bigr)$ is the list of principal stretches, or  singular values of $F\in\lp$. The following  is due to \cite{s}.

\begin{lemma}\label{l6}  [\cite{s}] A tensor $F_0\in\lp$ is a point of convexity of $W$ if and only if (i) $\la^0=\la(F_0)$ is a point of convexity of $\Phi$ and (ii)\beq \label{ss2}\frac{\partial\Phi}{\partial \la_i} (\la^0) + \frac{\partial\Phi}{\partial \la_j} (\la^0)\ge0,\quad i\neq j\in\{1,...,N\}\eeq
moreover in these circumstances, 
\beq\label{ord} \la_i^0\ge \la_j^0 \implies \frac{\partial\Phi}{\partial \la_i} (\la^0) \ge  \frac{\partial\Phi}{\partial \la_j}(\la^0)\eeq\end{lemma}

\begin{definition}\label{isp}
We call $\la\in\renp$ an \emph{isolated point of convexity of} $\Phi: \renp \to\re$ if it is a point of convexity of $\Phi$ and the latter is strictly convex in a neighborhood of $\la^0$.
\end{definition} 

\begin{proposition}\label{l7} Suppose that $N=2,3$, $W$ is isotropic, and that  (i) $\cD$ and $\cS$ satisfy the Partition Condition (Definition \ref{coo}), (ii) $W$ is strictly rank-one convex, (iii) $\la(F_0)$ is an isolated point of convexity of $\Phi$, (iv) $\Phi$ is pairwise nondecreasing, namely, inequalities \eqref{ss2} hold.
Then the only  classical equilibrium solution $y$ of Problem A is $y=F_0x$ on $\Om$.
\end{proposition}

\begin{remark}
Conditions (iii) and (iv) are weaker than in our previous results. Because we do not require the strict version of \eqref{ss2}, $F_0$ need not be a point of strict, or even isolated convexity of $W$, it is though a point of convexity.
\end{remark}
\bpr  {\bf Step 1.} We show  that any classical equilibrium solution $y$ has constant principal stretches, those of $F_0$, or $\la(\ny)=\la^0:=\la(F_0)$ on $\Om$. Let $\sio= {\partial\Phi} (\la^0) /{\partial \la_i} $
Define $\hw$ by \eqref{hw} and similarly let
$$\widehat\Phi (\la):=\Phi(\la)-\Phi(\la^0)-\sio(\la_i-\la_i^0) \quad \forall \la\in\renp$$
with summation over $i\in\{1,\ldots,N\}$. Now $F_0$ is a point of convexity of $W$ by virtue of Lemma \ref{l6} and Hypotheses (iii), (iv). Also $\la^0$ is an isolated point of convexity of $\widehat\Phi$, therefore
\beq\label{mi}\hw(F)\ge\hw(F_0)=0,\quad \ \widehat\Phi(\la)\ge \widehat\Phi(\la^0)=0 \quad \forall F\in \lp, \quad \forall\la\in\renp.\eeq
Because of \eqref{ord}, the $\sio$ are ordered as the $\la_i^0$. By  \eqref{ss2}, only  one of them can be  negative  and has the least absolute value. From \citet[Lemma 1.3]{ro1}, \citet[Proposition 18.3.2(2)]{s},
$$S(F_0)\cdot F\le\sio \la_i(F) \quad
\forall F\in\lp.$$
This and isotropy imply 
\beq\label{wf}\hw(F)\ge\widehat\Phi(\la(F))  \quad
\forall F\in\lp.\eeq
The first of \eqref{mi} allows an appeal to Lemma \ref{l3} which ensures that $\hw(\ny)=\hw(F_0)=0$
 on $\Om$. This, \eqref{wf} and the second of \eqref{mi} lead to $0=\hw(\ny)\ge\widehat\Phi(\ny)\ge0$, thus \beq\label{fo}\hw(\ny)=\widehat\Phi(\la(\ny))=0 \quad\hbox{ on } \Om.\eeq 
By Hypothesis (ii) and Lemma \ref{l3}, \eqref{ff0} holds.
By continuity of $x\mapsto\la(\ny(x))$ (with ordered singular values) the set $\la(\ny(\Om))$ is connected, consists of zeroes of $\widehat\Phi$ and  contains $\la^0$ by \eqref{ff0}, but the latter is an isolated zero of $\widehat\Phi$ by Hypothesis (iii). We conclude that
\beq\label{lao} \la(\ny(x))=\la^0\quad\forall x \in\Om.\eeq
Henceforth we write $\la_i$ and $\s_i$ in place of $\la_i^0$ and $\s_i^0$.

\noindent {\bf Step 2.} We show that $y$ is harmonic on $\Om$. Below we write $F$ in place of $\ny$. Let $\La=\hbox{diag}(\la)$ and $\Sigma= \hbox{diag}(\sigma)$, so that in view of \eqref{lao} and isotropy of $W$ there are two rotation fields $Q,R:\Om\to\ort$, such that
\beq\label{fs}F =Q\La R, \quad S(F)=Q\Sigma R .\eeq
Observe that 
\beq\label{sy}F^TS(F)=R^T \La\Sigma R\in \sym\eeq
Without loss of generality assume that $F_0=\La$ so that $S(F_0)=\Sigma$.
From \eqref{mi}, \eqref{fo} we infer that the derivative $\hs(F)$ of $\hw$ vanishes on $\Om$, hence 
\beq\label{so}S(F)=S(F_0)=\Sigma,\eeq
which combined with \eqref{sy} shows that $F^T\Sigma=\Sigma F$ on $\Om$. Since also $W(F)=0$ and the Eshelby tensor $P$ of \eqref{et} is solenoidal, we obtain
\beq\label{ds}0=\diver(F^TS(F))=\diver(F^T\Sigma)= \diver(\Sigma F)=\Sigma \diver F\eeq
\noindent {\bf Case 1.} Suppose $\s_i\neq0$. Since $\Sigma$ is then nonsingular this yields 
\beq\label{har} \diver F=\Delta y=0\quad\hbox{ on } \Om\eeq

\noindent {\bf Case 2 .} Suppose one of the $\s_i$ vanishes. Ordering them according to index we have $\s_3\ge\s_2>0=\s_1$ because of \eqref{ss2}, \eqref{ord}. Then $\Sigma\in\sym$ is positive semi-definite and has rank 2. From \eqref{fs}, \eqref{so}, $Q\Sigma R=\Sigma$ or
$$(QR)(R^T\Sigma R)=\Sigma$$
Since rank$\Sigma=2$, \cite[Lemma 2.1]{ro2} guaranties uniqueness of the proper orthogonal polar decomposition above hence $Q=R^T$ (and $R^T\Sigma R=\Sigma$) so that $F=R^T\La R\in \sym$ on $\Om$. Then $\diver F=\D$tr$F$ since 
$$(\diver F)_i=F_{ik,k}=F_{ki,k}=y_{k,ik}=y_{k,ki}=F_{kk,i}.$$
But since $F=R^T\La R$, tr$F$=tr$\Sigma=$const. on $\Om$ and \eqref{har} follows.

\noindent {\bf Case 3 .} If exactly two of the $\s_i$ vanish, then, ordering them according to index, \eqref{ss2}, \eqref{ord} dictate that 
$$\s_3>0=\s_2=\s_1.$$
Hypothesis (ii) implies the strict Baker-Ericksen Inequalities,
\beq\label{sbe} (\la_i\s_i-\la_j\s_j)(\la_i-\la_j)>0 \hbox{ if } \la_i\neq\la_j,\eeq 
(no sum) which  demand that $\la_1=\la_2$, since $\s_1=\s_2=0$. Hence,
\beq\label{sp}\La=\la_1 I+(\la_3-\la_1)E_0, \quad \Sigma=\s_3 E_0,\eeq
where $E_0=e\otimes e$ with $e$ a constant unit vector. Using \eqref{sy}, 
$$F^TS(F)=R^T\La\Sigma R=\la_3\s_3 E, \quad E=R^T E_0 R \hbox{ on } \Om.$$
The first equality in \eqref{ds} then shows that $\diver E=0$ on $\Om$. Hence the Cauchy Green tensor $$C=F^TF=R^T \La^2R= \la_1^2 I+(\la_3^2-\la_1^{2})E$$
from \eqref{fs}, \eqref{sp}, is divergence free. The identity
$$\diver C=\D|F|^2 +F^T \diver  F $$
and the fact that $|F|^2=$const. on $\Om$ from \eqref{lao} confirm \eqref{har}.

\noindent {\bf Case 4.} The remaining case is  $\s_i=0$ so that  $\Sigma=0$. It follows from \eqref{sbe} that all $\la_i$ are equal, thus $\Lambda=\la_1 I$ and by \eqref{fs}, $F=\la_1 Z$ for some $Z:\Om\to\ort$. Compatibility dictates that $Z=$const., hence in view of \eqref{ff0}, $F=F_0$ on $\Om$.

\noindent {\bf Step 3.} By \eqref{har}, $y\in C^\infty(\Om)$. Also \eqref{lao} implies that $|\ny|^2=$const. on $\Om$. Hence
$$0=\Delta\, |\ny|^2=(y_{k,m}y_{k,m})_{,jj}=2(y_{k,m}y_{k,mj})_{,j}$$
$$=2(y_{k,mj} y_{k,mj}+y_{k,m}y_{k,mjj})=2 y_{k,mj} y_{k,mj}+ 2y_{k,m}(y_{k,jj})_{,m}
$$ Now \eqref{har} reads $y_{k,jj} =0$ and thus $y_{k,mj} y_{k,mj} =|\D\D y|^2=0$ so that $\ny=$const. and equal to $F_0$ on $\Om$ in view of \eqref{ff0}.\epr

\section{The Mixed Displacement/Cauchy-Pressure Problem}\label{s5}

Here we assume isotropy; we choose $F_0$ to be a uniform dilatation and replace the dead-load conditions on $\cS$ by a Cauchy pressure condition. Define
\beq\label{ph}\ph(J):=W(J^{1/N}I)\quad\forall J>0\eeq
and let
$$\cC=\{F\vert F\in\lp, \; F=\delta R,\; \delta>0,\;R\in\ort\} $$
be the set of conformal deformation gradients,
so that $F\in\cC\implies W(F)=\ph({\rm det} F)$.   The following result, due to \cite{mi}, shows the importance of the function $\ph$ in \eqref{ph}, essentially the restriction of the stored energy to dilatations.
\begin{lemma}\label{l8} [\cite{mi}]. Given $J>0$ let $\cA_J=\{F|F\in\lp,  \det F=J\}$.
If $W$ is isotropic, strictly rank-one convex, and  $W(F)\to\infty$ as $|F|\to\infty$, then
$$W(F)>\ph(J) \quad\forall F\in \cA_J\backslash \cC.$$
\end{lemma} 

\begin{remark} Actually \cite{mi} also requires $W(F)\to\infty$ as $\det F\to 0$, however this is only  needed in case of non-strict rank-one convexity and non-smooth $W$. In fact all that is necessary in the above lemma is that for each $J>0$,  $W(F)\to\infty$ as $|F|\to\infty$  with $F\in \cA_J$. In particular, it in not necessary that $W(\delta I)\to\infty$ as $\delta\to\infty$ (dilatations)
\end{remark}
The Cauchy stress $T$ is such that $S(F)=T(F)\cof F$, with $\cof F=\det(F)F^{-T}$ for $F\in\lp$. It is not difficult to show that for a dilatation with $F_0=J_0^{1/N}I$, so that $\det F_0=J_0>0$, the corresponding Cauchy stress is
$$T_0=T(F_0)=\ph'(J_0)I,$$
a hydrostatic pressure or tension. Moreover, for a smooth deformation $y$ subject to Cauchy traction $T(\ny)n_*=T_0n_*$ on $\cS_*=y(\cS)$ (with $n_*$ the unit outward normal to $\cS_*$), the corresponding Piola traction on $\cS$ is 
 $$S(\ny)n=T_0 \cof(\ny)n .$$  
We now pose the mixed \emph{displacement/Cauchy Pressure problem}: 
\begin{problem2*}
Given $J_0>0$, find a classical equilibrium $y$ subject to the boundary conditions of displacement
\beq\label{bcup} y(x)= J_0^{1/N} x \quad \forall x \in \cD \eeq
and traction
\beq \label{bcsp} S(\ny)n=\ph'(J_0)\cof (\ny)n \quad \hbox{on } \cS \eeq
where $n$ is the unit outward normal on $\po$.
\end{problem2*}
Uniqueness for this problem is addressed in
\begin{proposition}\label{l9} Suppose that $W$ is isotropic, and (i) $\cD$ and $\cS$ satisfy the Partition Condition (Definition \ref{coo}), (ii) $W$ is strictly rank-one convex, (iii) $J_0$ is a point of convexity of $\ph$ from \eqref{ph}, that is,
$$\ph(J)\ge \ph(J_0)+\ph'(J_0)(J-J_0) \quad \forall J>0.$$
If $y$ is an equilibrium solution of Problem B, then $$y(x)= J_0^{1/N} x \quad \forall x \in \Om.$$
\end{proposition}
\bpr Define
$$ \hw(F):=W(F)-\ph'(J_0)(\hbox{det}F-J_0) \quad \forall F\in\lp,$$
and
$$\widehat\ph(J):=\hw(J^{1/N}I)= \ph(J)- \ph(J_0)-\ph'(J_0)(J-J_0) \quad \forall J>0.$$
Then $\hw$ is also isotropic and strictly rank-one convex and satisfies
$$\hw(F) \ge\hw({\rm det}F I)=\widehat\ph({\rm det}F )\ge\widehat\ph(J_0)=\hw(J_0^{1/N}I)=0 \quad\forall F\in\lp,$$
where we have used Lemma \ref{l8} and Hypotheses (ii) and (iii). Also in terms of the derivative of $\hw$,
$$\widehat S(F):=D\hw(F)=S(F)-\ph'(J_0)\cof F,$$
\eqref{bcsp} takes the form $\widehat S(\ny)n=0 $ on $\cS $. Consequently Lemma \ref{l3} with $W$ replaced by $\hw$ and $F_0= J_0^{1/N}I$ applies and ensures that for a classical equilibrium solution $y$ of Problem 2, $\hw(\ny)=\hw(F_0)=\widehat\ph(J_0)=0$ on $\Om$. 
But Lemma \ref{l8} demands that $\hw(\ny)>\widehat\ph({\rm det}\ny)$ unless $\ny\in\cC$. 
As a result $y$ is conformal. In addition, by \eqref{ff0}
 $\ny= J_0^{1/N}I$ on $\cD$. 
 Choose $z\in\cD$ with dist$(z,\cS)>\varepsilon>0$, let  $B=B_\varepsilon(z)$ (open ball) and extend $y$ onto $B\cup\Om$ 
  by letting  $y(x)= J_0^{1/N} x$ for $x\in B\backslash\Om$. 
 Then $y\in C^1(B\cup\Om)$ and conformal, hence  the result follows from the Identity Theorem. \epr

\section{Remarks and Examples}\label{s6}

\subsection{Implications}\label{60} In continuation or bifurcation studies, \cite{hs} for example, existence of solution branches to a parameter-dependent problem is investigated. The initial state is usually chosen to be the reference state (or another homogeneous deformation). Its unqualified  uniqueness can eliminate one of the Rabinowitz alternatives regarding boundedness of the  solution branch; see \cite[Theorem 4.2]{ks}.  Our uniqueness results do precisely that for the mixed problem, in case $\Om$ is diffeomorphic to a spherical shell with the Partition Condition in force, so that the outside surface is $\cD$ and the inside surface is $\cS$ (e.g., Fig.1~b). 

In some examples where the Partition Condition holds, $\cD$ comprizes two connected components (Fig.~1h,i).  Here our results would seem to be violated by rotating one of the components of $\cD$ by an angle of $2\pi$ about the vertical axis, while holding the other fixed (we suppose that $\Om$ is a solid of revolution about the vertical axis, Fig.~1h,i.) Actually, our results  imply nonexistence of  ``twisted'' classical solutions. If such equilibria exist, they must suffer a loss of classical smoothness, possibly with singularities developing at reentrant corners where $\cD$ and $\cS$ touch.

\subsection{Partition Condition Examples}\label{61} Our results apply to regions $\Om$ whose boundary satisfies the Partition Condition (Definition \ref{coo}). An obvious example is spherical shell, with $\cD$ the outside sphere (surface) and $\cS$ the inside one. More generally the difference of two  regions, $\Om=\Om_2\backslash\Om_1$ with $\Om_1\subset\Om_2\subset\ren$ and both star shaped with respect to the origin (Fig.~1a,b). Here $\cD=\po_2$ and $\cS=\po_1$. The restriction $\Om_1\subset\Om_2$ can be dropped however, see Fig.~1c. 
 
Given any  convex $\Om$ with the origin $O$ outside it, each ray from $O$ intersecting $\inter\Om$ will  intersect $\po$ at two points. Let the point closest to $O$ belong to $\cS$ and the farthest belong to $\cD$ (Fig.~1d,e). Otherwise if a ray touches $\po$ without intersecting the interior, the contact point can belong to either $\cS$ or $\cD$.
 
More generally, given \emph{any} (not necessarily star shaped) $\Om$ with piecewise smooth $\po$ and any choice of origin $\not\in\Om$, construct $\cD$ and $\cS$ as follows: if a ray from the origin enters 
$\Om$ at $x_0\in\po$ (except corners), let $x_0\in\cS$; if it leaves $\Om$ at $x_0\in\po$, let $x_0\in\cD$. If the ray is tangent to $\po$ at $x_0$, then the latter can belong to either subset (Fig.~1f)  As a result the Partition Condition does not restrict the region $\Om$ (which need not be star shaped), but rather the way in which $\po$ is partitioned into $\cD$ and $\cS$. 

Another possibility is that $\Om$ is the union of two regions star shaped with respect to the origin (Fig.~1h,i).

The Partition Condition does exclude some important configurations, such as a cylindrical region with $\cD$ consisting of the two opposite ends and $\cS$ the lateral boundary, of relevance in modeling the uniaxial test. On the other hand, even the slightest  tapering of the cylinder into a truncated cone, whose base is $\cD$ and the rest of $\po$ comprises $\cS$, ensures that the Partition Condition holds (Fig.~1g).

\subsection{Points of Convexity of the Stored Energy} \label{62}

While global convexity of $W$ is incompatible with the presence of a natural state and frame indifference (\citet{b}), we merely require in Propositions \ref{l5} and \ref{l7} that the homogeneous deformation gradient $F_0$  be a point of (sometimes strict) convexity of $W$ which need not be globally convex. In Lemma \ref{l3} and Proposition \ref{l4}, $F_0$ is assumed to be a natural state, which is a point of convexity. A point $F_0$ of strict convexity of $W$ as in Proposition \ref{l5} requires $\la^0$ to be a point of strict convexity of $\Phi$ and the strict version of \eqref{ss2}. For $N=2$ the latter allows states of uniaxial tension but excludes states of uniaxial compression. In particular, consider the compressible Neo-Hookean stored energy
$$W(F)= (1/2)|F|^2+(\det F)^{-\gamma}/\gamma \quad \forall F \in\lp$$
($N=2$, $\gamma\ge1$). This is strongly elliptic, hence strictly rank-one convex. One shows that every state $F_0$ with $\det F_0>1$  is a point of strict convexity of $W$, whereas states with $\det F_0<1$ are not even points of  convexity. This is because the Hessian matrix of $\Phi(\la_1,\la_2)=(\la_1^2+\la_2^2)/2+(\la_1\la_2)^{-\gamma}/\gamma$ is positive-definite on $\re_+^2$, whereas $\partial\Phi/\partial\la_1+\partial\Phi/\partial\la_2>0$ if and only if $\la_1\la_2>1$, as easily verified. This includes all states of uniaxial tension with $\partial\Phi/\partial\la_1=0$ and $\partial\Phi/\partial\la_2>0$, whereupon Proposition \ref{l5} applies. In contrast, uniaxial compression states ( $\partial\Phi/\partial\la_1=0$, $\partial\Phi/\partial\la_2<0$) are not covered by Proposition \ref{l5}, as they violate \eqref{ss2}, hence they are not points of convexity of $W$. It should be emphasized that this  is consistent with our expectations of non-uniqueness, due to the possibility of \emph{buckling} of the tapered column  in Fig.~1g under uniaxial compression loads in the vertical direction.  

For isotropic materials, strict convexity of $W$ at $F_0$ is slightly relaxed in Proposition \ref{l7}. Instead, we require that the singular value list $\la^0$ of $F_0$ be an \emph{isolated point of convexity} of $\Phi$ in \eqref{Phi} and  Definition \ref{isp}, and also that \eqref{ss2} hold. This can include uniaxial tension in three dimensions, but will always exclude uniaxial compression. The natural state of an isotropic $W$ is now allowed, although it is not a point of strict convexity of $W$ because of frame indifference.

Proposition \ref{l9} shows uniqueness for Problem B, where traction boundary conditions equivalent to Cauchy hydrostatic pressure or tension  are imposed on $\cS$ (Eq. \eqref{bcsp}). Hypothesis (iii) here is weaker than that of Proposition \ref{l7}. Suppose $J_0$ is a point of convexity  of $\ph$ (see \eqref{ph}) in the case of pressure, namely $\ph'(J_0)<0$. Then $F_0=J_0^{1/N}I$ is \emph{not} a point of convexity of $W$, as it  violates \eqref{ss2}, all principal stresses being equal and negative.  However, uniqueness still holds for Problem B though not necessarily for Problem A which prescribes the Piola traction due to the specified dilatation. Here \eqref{bcs} becomes $$S(\ny)n=J_0^{(N-1)/N}\ph'(J_0)n \hbox{ on } \cS ,$$
which is different from \eqref{bcsp}.

Convexity of $\ph$ does not follow from strict  rank-one convexity of $W$, or even strong ellipticity. An example with $\ph$ having the form of a two-well potential despite  $W$ being strictly polyconvex, is given by \citet[Eq. (6.6)]{rs}. This stored energy has two distinct dilatations that are natural states. They are the global minima, hence points of convexity of $\ph$ of \eqref{ph}. As a result, Proposition \ref{l9} still applies to $F_0$ chosen to be one of them, with null traction specified on $\cS$.
 
 \subsection{Partition Condition Omission and Nonuniqueness}\label{63}
We provide an example of nonuniqueness in case the Partition Condition is omitted from our hypotheses. It is well known that the pure dead-load traction problem suffers from lack of uniqueness, so we consider a genuine mixed problem, where both $\cD$ and $\cS$ have nonempty relative interior.

Let $N=2$ and  $\Om$ be the annulus $\{x|x\in\rt, \;\; 1\le x\le 5\}$. Suppose $\cD$ is the inside circle $|x|=1$ and $\cS$ the exterior one $|x|=5$, whereby both inequalities of the Partition Condition are reversed (violated). Choose
\beq\label{wsp}W(F)=\sqrt{|F|^2+2\det F} +\frac{k}{2} \Bigl(\det F-\frac{1+k}{k}\Bigr)^2 \quad \forall F\in \lp,\eeq
with $k>0$ a constant.
This energy is isotropic; it is strongly elliptic as it satisfies 
\cite[Eq. (1.35)]{jkkes}, 
hence strictly rank-one convex, and polyconvex by 
\cite[Proposition 3.1]{rs}. 
Moreover, for $k\ge 2$ it is globally minimized on $\lp$ at $\ort$, while $S(I)=0$ for any $k>0$.  We impose $ y(x)=x$ for $|x|=1$ and $S(\ny)x/5=0$ for $|x|=5$. This is consistent with null Piola traction conditions \eqref{bcs} for  Problem A and also null Cauchy pressure conditions \eqref{bcsp} for Problem B. The natural state $F_0=I$  minimizes $W$ together with all rotations provided $k\ge 2$. Also $(1,1)$ is an isolated point of convexity of $$\Phi(\la_1,\la_2)=W(\rm{diag}(\la_1,\la_2))=\la_1+\la_2+ \frac{k}{2} \Bigl(\la_1\la_2-\frac{1+k}{k}\Bigr)^2, \quad \la_i>0, \;i\in\{1,2\} $$ (Definition \ref{isp}). In addition, the dilatation energy function \eqref{ph} now becomes $$\ph(J)=2\sqrt{J}+ \frac{k}{2} \Bigl(J-\frac{1+k}{k}\Bigr)^2 .$$ For $k\ge 2$, $\ph(J)$ has a global minimum at $J=1$, which is therefore a point of its convexity (although $\ph$ is not globally convex). As a result, $F_0=I$ with $J_0=1$ satisfies Hypothesis (iii) of Propositions \ref{l6}, \ref{l7} and \ref{l9}.  One solution is clearly the identity map $y(x)=x$ for all $x\in\Om$. 

We construct another radial solution of the form $y(x) = R(r)x/r$, where $r=|x|\in[1,5]$, the homogeneous one corresponding to $R(r)=r$. Since the principal stretches (singular values of $\ny$) are 
\beq\label{lir}\la_1(r)=R'(r), \quad\la_2(r)=R(r)/r,\eeq  let 
$$\Phi_i(r)=\frac{\partial\Phi}{\partial \la_i}\bigl(R'(r),R(r)/r\bigr).$$ 
As is well known, $y$  is a radial solution of the problem if and only if $R:[1,5]\to\re$ satisfies the radial equilibrium ODE 
\beq\label{ode}r\Phi'_1(r)+\Phi_1(r)-\Phi_2(r)=0, \quad 1<r<5\eeq 
and the boundary conditions
\beq\label{rbc}R(1)=1,\quad \Phi_1(5)=0.\eeq
In fact (\cite{hm}), a radial solution of \eqref{ode}  is 
\beq\label{solr} R(r)=\sqrt{Jr^2+b},\eeq
where $J$ and $b$ are constants with $J>0$ and $J+b>0$.   It follows that $\det\ny=J=$const. on $\Om$. This solution is universal to the class of Varga materials of the form 
$$\Phi(\la_1,\la_2)=\la_1+\la_2+ \psi\Bigl(\la_1\la_2\Bigr), \quad \la_i>0, \;i\in\{1,2\} $$ 
(\citet{hm}) for a smooth $\psi:\re_+\to\re$ with $\psi'(1)=-1$ for $F_0=1$ to be stress free. Here  $\psi(J)= \frac{k}{2} \Bigl(J-\frac{1+k}{k}\Bigr)^2$. The boundary conditions \eqref{rbc} become 
$$R(1)=1, \quad 1+\psi'(J)R(5)/5=0$$
where we have used \eqref{lir}. By \eqref{solr} they reduce to
\beq\label{b}b=1-J,\eeq and $J>0$  a root of
\beq\label{root} \psi'(J)=f(J), \quad \hbox{ where } f(J):=-\sqrt{\frac{25}{1+24J}}.\eeq
 Note that $J=1$ is always a root of \eqref{root},  so that $b=0$ by \eqref{b}, and \eqref{solr} reduces to the homogeneous solution $R(r)=r$ or $y=x$. Here $\psi'(J)= k \Bigl(J-\frac{1+k}{k}\Bigr)$ is linear and strictly increasing. The right-hand side  $f(J)$ of \eqref{root} is negative, strictly increasing from it value $f(0)=-5$, concave, and its graph intersects that of $\psi'$ at $(1,-1)$, with slope $f'(1)=12/25$. Thus for any value of $k=\psi''$ such that $12/25<k<4$ the straight graph of $\psi'$ will intersect that of $f$ at  second root $J=J_*\in(0,1)$, with 
 $$k=\frac{-1-f(J_*)}{1-J_*}.$$ We require $k\ge 2$ for $F_0$ to be a natural state.  For example, $R(r)=\sqrt{r^2/20+19/20}\;$ for $1\le r \le 5$ is a solution with $J_*=1/20$ for $k\approx 2.496$. This shows nonuniqueness for Problems A and B when the 
 Partition Condition is violated, but all other hypotheses of Propositions \ref{l6}, \ref{l7} and \ref{l9} hold. This inhomogeneous deformation is not a global minimizer of the total energy, unlike the homogeneous one.
 
We remark that when we interchange $\cD$ and $\cS$ for the annulus, and prescribe the identity map at $r=5$ and null tractions at $r=1$, the Partition Condition holds. Not surprisingly, we now find no radial solutions in the previous example other than the identity $R(r)=r$, because the analogue of the function $f$ of \eqref{root} is  strictly decreasing and can only have one intersection with the increasing graph of $\psi'$, which occurs at $J=1$.
 
 In this case, another way to show uniqueness of radial deformations is that \eqref{ff0} which now holds, together with the boundary condition  $R(5)=5$ determines $R'(5)=1$ as well, thus these  Cauchy initial conditions  for \eqref{ode} only allow the homogeneous solution. See \cite{ks}. Of course our results also preclude non-radial solutions as well.

\section*{Acknowledgements}
This work is dedicated to Rohan Abeyaratne, scholar, teacher, academic brother and dear friend, on the occasion of his 
$70^{\rm th}$ birthday, in recognition of his seminal contributions to the understaning of Nonlinear Elasticity. 

\bibliographystyle{model1-num-names}
\bibliography{Refsuniq}

\begin{thebibliography}{12}
\expandafter\ifx\csname natexlab\endcsname\relax\def\natexlab#1{#1}\fi
\providecommand{\bibinfo}[2]{#2}
\ifx\xfnm\relax \def\xfnm[#1]{\unskip,\space#1}\fi
%Type = Article
\bibitem[{Knops and Stuart(1984)}]{ks}
\bibinfo{author}{R.~Knops}, \bibinfo{author}{C.~A. Stuart},
\newblock \bibinfo{title}{Quasiconvexity and uniqueness of equilibrium solutions in nonlinear elasticity},
\newblock \bibinfo{journal}{Archive for rational mechanics and analysis} \bibinfo{volume}{86} (\bibinfo{year}{1984}) \bibinfo{pages}{233--249}.
%Type = Article
\bibitem[{Sivaloganathan and Spector(2018)}]{sis}
\bibinfo{author}{J.~Sivaloganathan}, \bibinfo{author}{S.~J. Spector},
\newblock \bibinfo{title}{On the uniqueness of energy minimizers in finite elasticity},
\newblock \bibinfo{journal}{Journal of Elasticity} \bibinfo{volume}{133} (\bibinfo{year}{2018}) \bibinfo{pages}{73--103}.
%Type = Article
\bibitem[{Ball(1976)}]{b}
\bibinfo{author}{J.~M. Ball},
\newblock \bibinfo{title}{Convexity conditions and existence theorems in nonlinear elasticity},
\newblock \bibinfo{journal}{Archive for rational mechanics and Analysis} \bibinfo{volume}{63} (\bibinfo{year}{1976}) \bibinfo{pages}{337--403}.
%Type = Book
\bibitem[{Silhavy(2013)}]{s}
\bibinfo{author}{M.~Silhavy}, \bibinfo{title}{The mechanics and thermodynamics of continuous media}, \bibinfo{publisher}{Springer Science \& Business Media}, \bibinfo{year}{2013}.
%Type = Article
\bibitem[{Ball and James(1987)}]{bj}
\bibinfo{author}{J.~M. Ball}, \bibinfo{author}{R.~D. James},
\newblock \bibinfo{title}{Fine phase mixtures as minimizers of energy},
\newblock \bibinfo{journal}{Archive for Rational Mechanics and Analysis} \bibinfo{volume}{100} (\bibinfo{year}{1987}) \bibinfo{pages}{13--52}.
%Type = Article
\bibitem[{Rosakis(1997)}]{ro1}
\bibinfo{author}{P.~Rosakis},
\newblock \bibinfo{title}{Characterization of convex isotropic functions},
\newblock \bibinfo{journal}{Journal of elasticity} \bibinfo{volume}{49} (\bibinfo{year}{1997}) \bibinfo{pages}{257--267}.
%Type = Article
\bibitem[{Rosakis(1990)}]{ro2}
\bibinfo{author}{P.~Rosakis},
\newblock \bibinfo{title}{Ellipticity and deformations with discontinuous gradients in finite elastostatics},
\newblock \bibinfo{journal}{Archive for Rational Mechanics and Analysis} \bibinfo{volume}{109} (\bibinfo{year}{1990}) \bibinfo{pages}{1--37}.
%Type = Article
\bibitem[{Mizel(1998)}]{mi}
\bibinfo{author}{V.~J. Mizel},
\newblock \bibinfo{title}{On the ubiquity of fracture in nonlinear elasticity},
\newblock \bibinfo{journal}{Journal of elasticity} \bibinfo{volume}{52} (\bibinfo{year}{1998}) \bibinfo{pages}{257--266}.
%Type = Article
\bibitem[{Healey and Simpson(1998)}]{hs}
\bibinfo{author}{T.~J. Healey}, \bibinfo{author}{H.~C. Simpson},
\newblock \bibinfo{title}{Global continuation in nonlinear elasticity},
\newblock \bibinfo{journal}{Archive for rational mechanics and analysis} \bibinfo{volume}{143} (\bibinfo{year}{1998}) \bibinfo{pages}{1--28}.
%Type = Article
\bibitem[{Rosakis and Simpson(1994)}]{rs}
\bibinfo{author}{P.~Rosakis}, \bibinfo{author}{H.~C. Simpson},
\newblock \bibinfo{title}{On the relation between polyconvexity and rank-one convexity in nonlinear elasticity},
\newblock \bibinfo{journal}{Journal of elasticity} \bibinfo{volume}{37} (\bibinfo{year}{1994}) \bibinfo{pages}{113--137}.
%Type = Article
\bibitem[{Knowles and Sternberg(1978)}]{jkkes}
\bibinfo{author}{J.~K. Knowles}, \bibinfo{author}{E.~Sternberg},
\newblock \bibinfo{title}{On the failure of ellipticity and the emergence of discontinuous deformation gradients in plane finite elastostatics},
\newblock \bibinfo{journal}{Journal of Elasticity} \bibinfo{volume}{8} (\bibinfo{year}{1978}) \bibinfo{pages}{329--379}.
%Type = Article
\bibitem[{Horgan and Murphy(2009)}]{hm}
\bibinfo{author}{C.~O. Horgan}, \bibinfo{author}{J.~G. Murphy},
\newblock \bibinfo{title}{Some analytic solutions for plane strain deformations of compressible isotropic nonlinearly elastic materials},
\newblock \bibinfo{journal}{Advances in Mathematical Modeling and Experimental Methods for Materials and Structures: The Jacob Aboudi Volume}  (\bibinfo{year}{2009}) \bibinfo{pages}{237--247}.

\end{thebibliography}

%\appendixname ejchbducbducuee efecdcdecd
\end{document}